\documentstyle[12pt]{article}
\newlength{\abstractwidth}
\renewcommand{\thefootnote}{\fnsymbol{footnote}}
\renewcommand{\thanks}[1]{\footnote{#1}} 
\newcommand{\starttext}{
\setcounter{footnote}{0}
\renewcommand{\thefootnote}{\arabic{footnote}}}

\newcommand{\be}{\begin{equation}}
\newcommand{\bea}{\begin{eqnarray}}
\newcommand{\eea}{\end{eqnarray}}
\newcommand{\ee}{\end{equation}}

\def\ba{\begin{eqnarray}}
\def\ea{\end{eqnarray}}

\def\D{{\cal D}}

\def\K{{\cal K}}

\def\z{{\bf z}}

\def\log{\,{\rm log}\,}

\def\o{\omega}

\def\al{\alpha}
\def\b{\beta}
\def\g{\gamma}
\def\d{\delta}
\def\e{\epsilon}

\def\l{\lambda}
\def\m{\mu}

\def\o{\omega}

\def\r{\rho}

\def\z{\zeta}
\def\D{\Delta}
\def\O{\Omega}

\def\R{{\bf R}}
\def\C{{\bf C}}
\def\P{{\bf P}}

\def\i{\infty}
\def\I{\int}

\def\s{\sum}

\def\ddb{{\partial\bar\partial}}

\def\sub{\subseteq}
\def\ra{\rightarrow}

\def\D{\Delta}

\def\cH{{\cal H}}

\def\cH{{\cal H}}

\def\us{{\underline s}}
\def\K{{K\"ahler\ }}

\def\v{\vskip .1in}

\def\[{{\bf [}}
\def\]{{\bf ]}}

\def\pl{\partial}

\begin{document}
\starttext
\baselineskip=18pt
\setcounter{footnote}{0}
\newtheorem{theorem}{Theorem}
\newtheorem{lemma}{Lemma}
\newtheorem{definition}{Definition}
\begin{center}
{\Large \bf THE MONGE-AMP\`ERE OPERATOR AND
GEODESICS IN THE SPACE OF
K\"AHLER POTENTIALS}\footnote{Research supported in part by
National Science Foundation
grants  DMS-02-45371 and DMS-01-00410}
\\
\bigskip

{\large D.H. Phong$^*$ and
Jacob Sturm$^{\dagger}$} \\

\bigskip

$^*$ Department of Mathematics\\
Columbia University, New York, NY 10027\\

\v

$^{\dagger}$ Department of Mathematics \\
Rutgers University, Newark, NJ 07102

\end{center}

\baselineskip=15pt
\setcounter{equation}{0}
\setcounter{footnote}{0}

\section{Introduction}
\setcounter{equation}{0}
Let  $X $ be an $n$-dimensional compact complex
manifold,
$L\ra X$ a positive holomorphic line bundle,
and $\cH$ the space of positively curved hermitian
metrics on $L$.
The purpose of this article is to prove that
geodesics in the infinite-dimensional symmetric
space
$\cH$ can be
uniformly approximated by geodesics in the
finite-dimensional
symmetric spaces
$\cH_k=GL(N_{k}+1,\C)/U(N_k+1,\C)$,
where $N_k+1=\dim(H^0(L^k))$.  Thus the
$\cH_k\sub\cH$ are becoming flat as $k\to\i$.

\v

The motivation for this work
comes from Donaldson's far reaching program
\cite{D97, D99} relating the geometry of ${\cal H}$
to the existence and uniqueness of constant scalar
curvature K\"ahler metrics. As advocated by S.T. Yau
over the years, ${\cal H}$ should be approximated by
${\cal H}_k$, and the properties of this approximation
should be closely reflected in many basic questions of 
K\"ahler geometry.
In particular, the condition of
``stability" is one which concerns
the growth of  energy functionals along
the geodesics of $\cH_k$. On the
other hand, the existence and uniqueness of metrics
of constant scalar curvature concerns
the growth of energy functionals along
the geodesics of $\cH$ (see \cite{D01,D02,T97,Y78,
Y}). Thus a good
understanding of
the relationship between these different
types of geodesics is desirable.

\v
In order to state the precise theorem, we need
some notation: Let $h:L\ra [0,\i)$ be a smooth
hermitian metric. If $s\in L$ we
write $h(s)=|s|_h$ and for $k>0$, we denote by
$h^k$ the induced metric on $L^k$. The
   curvature of $h$ is the $(1,1)$ form on
$X$  defined locally by the formula
$R(h)=-{\sqrt{-1}\over 2}\ddb \log|s(z)|_{h}^2$,
where $s(z)$ is a local,
nowhere vanishing holomorphic section.
In particular,
$R(h^k)=kR(h)$. Let
$$ \cH\ = \ \{h:L\ra [0,\i) : h
\hbox{\ \ is  smooth metric
on $L$ such that $R(h)>0 $\}}\ .
$$
If we fix $h_0\in\cH$ and let $\o_0=R(h_0)$ we have
a natural isomorphism
$$ \cH\ = \ \{\phi\in C^\i(X) : \o_\phi\ = \
\o_0+{\sqrt{-1}\over 2}\ddb\phi>0\ \} \ ,
$$
where $\phi$ is identified with $h=h_0e^{-\phi}$
so that $R(h)=\o_\phi$.
Then $\cH$ is an infinite dimensional manifold
whose tangent space $T_\phi\cH$ at
$\phi\in\cH$ is  naturally identified with
$C^\i(X)$.
\v
Now let $\phi\in\cH$ and $\psi\in C^\i(X)$ and
define a metric on $\cH$ by
\be\label{metric}
||\psi||^2_\phi\ = \ \I_X |\psi|^2\
\o_\phi^n\ \ .
\ee
Donaldson \cite{D99}, Mabuchi \cite{M87} and Semmes
\cite{S92}
have shown
that (\ref{metric}) defines a Riemannian
metric which makes $\cH$ into an infinite
dimensional negatively curved symmetric space.
Furthermore,
the geodesics of $\cH$ in this metric are the paths $\phi_t$ which 
satisfy the
partial differential equation
\be
\ddot\phi-|\partial\dot\phi|_{\o_\phi}^2=0.
\ee

\bigskip

The space $\cH$ contains a canonical family of
finite-dimensional negatively curved symmetric
spaces
$\cH_k$ which  are defined as follows: For $k>>0$
and for
$\us=(s_0,....,s_{N_k})$ an ordered basis of
$H^0(L^k)$, let
$$ \iota_\us: X\hookrightarrow \P^N
$$
be the Kodaira embedding given by
$z\mapsto (s_0(z),..., s_{N_k}(z))$. Then we have a
canonical isomorphism
$L^k=\iota_\us^*O(1)$. Let $h_{FS}$
be the Fubini-Study metric on $O(1)\ra \P^N$ and
let
\be\label{fubini}
h_\us=
(\iota_\us^*h_{FS})^{1/k}\
=
\  {h_0\over
\left(\s_{j=0}^{N_k} |s_j|^2_{h_0^k}\right)^{1/k}}
\ee
Note that the right side of (\ref{fubini}) is
independent of the choice of $h_0$. In
particular
\be\label{fubini two}
\s_{j=0}^{N_k} |s_j|^2_{h_\us^k}\ = \ 1
\ee
Let
$$ \cH_k\ = \ \{h_\us: \us \hbox{\ a basis
of  $H^0(L^k) $\} }\ \sub \ \cH
$$

\v
Then
$\cH_k = GL(N_k+1)/U(N_k+1)$ is a finite-dimensional 
negatively curved symmetric space
sitting inside of $\cH$.  We note that $\cH_k\sub \cH_l$ if $k|l$ and therefore,
if $k_1,k_2>0$ there exists $l$ such that
$\cH_{k_1}\sub \cH_l$ and
$\cH_{k_2}\sub \cH_l$. It
is well known that the $\cH_k$ are  topologically
dense in $\cH$: If $h\in \cH$ then there
exists
$h(k)\in\cH_k$ such that $h(k)\ra h$ in the
$C^\i$ topology. This follows from the Tian-Yau-Zelditch
theorem  on
the
density of states (Yau\cite{Y87}, Tian \cite{T90} and
Zelditch\cite{Z}; see also Catlin
\cite{C} for corresponding
results for the Bergman kernel). In fact, if
$h\in\cH$, then
there is a canonical choice of the approximating
sequence $h(k)$:
Let $\us$ be a basis
of $H^0(L^k)$ which is orthonormal with respect
to the metrics $h$.  In
other words,
\be\I_X (s_i,s_j)_{h^k} \
\o^n=\d_{ij}\ \ \hbox{where $\o=R(h) $} \ .
\ee
Note that $\underline s$ is not unique: if $\underline s$ is
orthonormal, so is $u\,\underline s$, for any $u\in U(N_k+1)$.
Define $\r_k(h)=\r_k(\o) =\s_j |s_j|^2_{h^k}$.
Then  Theorem 1 of \cite{Z}, which is the $C^\i$ version
of the $C^2$ approximation result first estsablished
in  \cite{T90}, says
that for
$h$ fixed, we have a $C^\i$ asymptotic expansion as
$k\to\i$:
\be\label{zelditch} \r_k(\o)\ \sim \ k^n +
A_1(\o)k^{n-1} + A_2(\o)k^{n-1}+\cdots
\ee
Here the $A_j(\o)$ are smooth functions on $X$
defined locally by $\o$ which can be computed in terms
of the curvature of $\o$ by the work of Lu \cite{Lu}. Let
$\hat\us=k^{-n/2}\us$ and $h(k)=h_{\hat\us}$. In particular,
(\ref{fubini}) and (\ref{zelditch}) imply that for each $r>0$
\be\label{powers of k}
\left|\left|{h(k)\over
h}-1\right|\right|=O\left({1\over k^2}\right)\ \ , \ \
||\o(k)-\o ||= O\left({1\over
k^2}\right)\ \ , \ \
||\phi(k)-\phi ||= O\left({1\over
k^2}\right)
\ee
where the norms are all taken with respect to
${C^r(\o_0)}$. Here, as before, $\o=R(h)$,
$\o(k)=R(h(k))$, $h=h_0e^{-\phi} $ and
$h(k)=h_0e^{-\phi(k)}$.

\bigskip

Now let $h_0,h_1\in\cH$.  It is known by the
work of Chen \cite{Ch} (see also
more recent progress in Donaldson \cite{D04}
and Chen-Tian \cite{CT}) that there is a unique
$C^{1,1}$ geodesic $h:[0,1]\ra \cH$ joining
$h_0$ to $h_1$. As discussed in Donaldson \cite{D99} 
the optimal regularity properties of this geodesic
are of considerable interest. Our main
theorem says that
$h_t=h_0e^{-\phi(t)}$ is a uniform limit of
the smooth geodesics in $\cH_k=GL(N_k+1)/U(N_k+1)$
which join $h_0(k)=h_{\hat \us^{(0)}}$ to
$h_1(k)=h_{\hat
\us^{(1)}}$.

\v
To be precise,  let
$\sigma\in GL(N_k+1)$ be the change of basis
matrix defined by
$\sigma\cdot\hat\us^{(0)}=\hat\us^{(1)}$. Without loss of
generality, we may assume that
$\sigma$ is diagonal with entries
$e^{\l_0},...,e^{\l_N}$ for some
$\l_j\in\R$. Let $\hat\us^{(t)}=\sigma^t\cdot\hat\us^{(0)}$
where $\sigma^t$ is diagonal with entries
$e^{\l_jt}$. Define
\be
h(t;k)=h_{\hat\us^{(t)}}=h_0e^{-\phi(t;k)}.
\ee
Then $h(t;k)$ is the smooth geodesic in
$GL(N_k+1)/U(N_k+1)$ joining $h_0(k)$ to
$h_1(k)$. Explicitly, using
(\ref{fubini}), we can also write
\be
\label{phi}
   \phi(t;k)\ = \ {1\over k}
\log\left(\s_{j=0}^N e^{2\l_jt}|\hat s_j^{(0)}|^2_{h_0^k}
\right).
\ee
 
\v
\begin{theorem} Let $h_0,h_1\in \cH$ and
$h_t=h_0e^{-\phi_t}$  the unique $C^{1,1}$
geodesic joining $h_0$ to $h_1$. Then
\be
\phi_t\ = \ \lim_{l\to \i}
[\sup_{k\geq l} \phi(t;k)]^* \ \
\hbox{uniformly as $l\to\i $},
\ee
where, for any bounded function
$u: X\times [0,1]\ra \R$, we define
the upper envelope $u^*$ of $u$ by
$u^*(\z_0)=\lim_{\e\to 0}\sup_{|\z-\z_0|<\e} u(\z)
$.
\end{theorem}
\v
{\it Remark 1}. The proof will show that
\be
\phi_t\ = \ \lim_{l\to \i}
[\sup_{k\geq l} \phi(t;k)]\ \ \ \hbox{almost
everywhere.}
\ee

\v
{\it Remark 2}. We note that
the upper envelope $u^*$ is independent of the choice
of coordinate systems defining the balls
$|\z-\z_0|<\e$. It is the smallest
upper semicontinuous function which is greater than
or equal to $u$. If $u_k$ is a sequence of plurisubharmonic
functions which are locally uniformly bounded, then
$[{\rm sup}\,u_k]^*$ is plurisubharmonic and equal
to ${\rm sup}\,u_k$ almost everywhere.
Similarly, we also define for later use the lower envelope
$u_*$ of a bounded function $u$ by
$u^*(\z_0)=\lim_{\e\to 0}\inf_{|\z-\z_0|<\e} u(\z)$.
The function $u_*$ is the largest lower semi-continuous
function which is less than or equal to $u$.
\v

To prove the theorem we first apply the observation
of Donaldson \cite{D99}, Mabuchi \cite{M87} and Semmes
\cite{S92}, which shows that solving the geodesic
equation on $\cH$ is equivalent to solving the
degenerate Monge-Amp\`ere equation
\be
\label{MA}
\O_\Phi^{n+1} = 0 \ \
\ee
on the manifold $\bar M=X\times A$, where $A$ is
the annulus $A=\{w\in\C: 1\leq |w|\leq e\}$,
the values of $\Phi$ are prescribed on the
boundary of $\bar M$ by smooth rotationally
symmetric data, and $\Omega_\Phi|_{X\times\{t\}}$
is positive for every $t$. Here
we are writing
$\O_\Phi=\O_0+{\sqrt{-1}\over 2}\ddb\Phi$ where
$\O_0=\pi_1^*\o_0$ and $\Phi$ is a smooth on $\bar
M$.  We then use the Tian-Yau-Zelditch theorem to prove that
the $\cH_k$ geodesics are, in a certain sense,
approximate solutions to (\ref{MA}). A key
step is then proving that the limit of
the sequence of approximate solutions is
in fact a weak solution of the Monge-Amp\`ere equation
and that the weak solution thus obtained is unique. This
is accomplished using  the  methods of pluripotential
theory which were introduced and systematically developed
by Bedford-Taylor in their fundamental
work \cite{BT76} and \cite{BT82}. In recent years, this
subject has been extended by  Demailly \cite{D},
Klimek \cite{K}, Blocki \cite{B}, and others. The uniqueness
implies that the limit coincides with
the $C^{1,1}$ solution $\phi_t$. Now the uniqueness
theorem for bounded open sets in $\C^n$, which is due to
Bedford-Taylor, requires smoothing techniques for rough
solutions which do not generalize in a simple way to
manifolds. It is quite possible that the 
regularization methods of Demailly, which have
recently been successfully applied by Guedj and
Zeriahi \cite{GZ} in the setting of compact manifolds,
will also apply in our setting. However, we use a
different approach which exploits the  particular
structure of the weak solutions in our case.
A particularly important ingredient for us is the existence of
a vector field $Y$, transversal to the boundary of $\bar M$,
with
\be
|Y(\Phi_k)|\leq C,
\ee
where $\Phi_k$ is a sequence of approximate solutions to
the Dirichlet problem for the Monge-Amp\`ere equation
(c.f. Theorem 3 and the proof of Theorem 6
below). For the application to Theorem 1,
this hypothesis is a consequence of the a priori estimate
\be
|\dot\phi(t,k)|\leq C,
\ee
which is the one global estimate for the derivatives of $\phi(t;k)$
that we can actually obtain. 
\v
One approach to the question raised in Donaldson 
 \cite{D99}  on
the optimal regularity of geodesics is the establishment
of a priori estimates on the smooth approximate
solutions $\phi(t;k)$. As we just noted,
the uniform $C^0$ estimates for $\phi(t;k)$
and for $\dot\phi(t;k)$ do hold.
But extending these to higher derivative
estimates seems rather hard, and at the
same time, quite intriguing: For example, the
estimate which is needed for
$\ddot\phi(t;k)$ can be carried
out in certain special cases, and appears to be
related to the central limit theorem in
probability theory. In the last section, we
shall make some  remarks along these
lines.

\section{The volume estimate}
\setcounter{equation}{0}

In this section, we establish the basic properties of the functions
$\phi(t;k)$ introduced in (\ref{phi}).
The properties of the functions $\phi(t;k)$ themselves
are summarized in Lemma 1 below, while the properties
of their corresponding Monge-Amp\`ere measures are
given in Theorem 2.

\bigskip

$\bullet$ A first key ingredient is the following very
general volume formula.
Let $(X,\o_0)$ be a compact \K manifold of dimension $n$, 
$A=\{w\in \C: 1\leq |w|\leq e\}$,
$\bar M=X\times A$,
$\pi_1:\bar M\to X$ be the projection on the
first factor, and $\O_0=\pi_1^*\o_0 $. Then $\O_0$ is a
closed positive
$(1,1)$-form on $\bar M$ such that $\O_0^{n+1}=0$.

\v
Let $\phi:[0,1]\ra \phi(t)\in\cH
=\{\phi\in C^\i(X): \o_0+{\sqrt{-1}\over 2}\ddb\phi\}$
be a smooth path joining $\phi(0)$
to $\phi(1)$ and define $\Phi: \bar M\ra \R $ by
\be\label{Phi}  \Phi(z,w)\ = \ \phi(t)(z)\ \ \hbox{
\ where $t=\log|w|$
}
\ee
and the corresponding $(1,1)$-form $\Omega_\Phi$ by
\be
\label{OPhi}
\O_\Phi\ = \ \O_0+{\sqrt{-1}\over 2}\ddb \Phi.
\ee
In local coordinates $z^i$ for $X$ and $v=\log w$ for $A$,
if we identify a K\"ahler form $\o_0={\sqrt{-1}\over 2}\sum_{ij}g_{\bar 
j i}^0dz^i\wedge d\bar z^j$
with the hermitian matrix $\{g_{\bar j i}^0\}$, then we can write
$$ \O_\Phi\ = \ \left(
\matrix{g_{\bar ji}^0+\pl_i\pl_{\bar j}\phi & {1\over 2}B\cr
            {1\over 2}B^*& {1\over 4}\ddot\phi\cr
}
\right)
$$
where $B$ is the column vector
$$ B\ = \ \left(
\matrix{{\pl \dot \phi\over \pl \bar z_1}\cr
                {\pl  \dot\phi\over \pl \bar z_2}\cr
                \cdot\cr
\cdot\cr
{\pl \dot \phi\over \pl \bar z_n}\cr
            }
\right)
$$
and $B^*$ is the conjugate transpose of $B$. It follows that
\be 
\O_\Phi^{n+1}\ = \ {1\over 4}(\ddot\phi-|\pl\dot\phi|^2_{\o_\phi})\
\o_\phi^n\wedge ({\sqrt{-1}\over 2}dv\wedge d\bar v),
\ee
and the condition $\O_\Phi^{n+1}=0$ is equivalent to
the geodesic equation $\ddot\phi-|\pl\dot\phi|^2_{\o_\phi}=0$.
This is a key observation which was pointed out in
 \cite{D99}, \cite{M87} and 
\cite{S92}.
Now
\be
\label{volume0}
\int_{X\times A}\O_\Phi^{n+1}\ = \
\I_0^1\I_X(\ddot\phi-|\pl\dot\phi|^2)\ \o_\phi^n\,dt\ = \ \I_0^1
{d\over dt}\bigg(\int_X\dot\phi\o_\phi^n\bigg)\ dt.
\ee
and we obtain the desired volume formula:
\be
\label{volume}
\int_{X\times A}\O_\Phi^{n+1}\ = \
\I_X \dot\phi(1)\o_{\phi(1)}^n\ - \
\I_X \dot\phi(0)\o_{\phi(0)}^n.\
\ee
Note that  (\ref{volume}) is only a true volume if
$\O_\Phi\geq 0$, but that (\ref{volume}) is valid
even in the absence of this hypothesis.

\v
\bigskip
$\bullet$ It is perhaps noteworthy that ${1\over 
V}\int_X\dot\phi\o_\phi^n=\dot
F_{\o_{\phi_0}}^0$, where $V=\int_X\o_0^n$,
and $F_{\o_{\phi_0}}^0(\phi)$
is precisely the functional whose critical points
(with the constraint ${1\over V}\int_Xe^{f_0-\phi}\o_0^n=1$,
and $f_0$ is a fixed function satisfying
$Ric(\o_0)-\mu\o_0={\sqrt{-1}\over 2}\ddb f_0$),
give K\"ahler-Einstein metrics. Thus the preceding relation can be
rewritten as
\be
\label{F0}
{1\over V}\I_{X\times A}\O_\Phi^{n+1}\ = \
\int_0^1\ddot F_{\o_{\phi_0}}^0dt
\
=
\ \dot F_{\o_{\phi_0}}^0(\phi(1))
-
\dot F_{\o_{\phi_0}}^0(\phi(1)).
\ee
The asymptotic behavior
of $F_{\o_{\phi_0}}^0$ along geodesics in ${\cal H}_k$
is known to be closely related to Chow-Mumford stability
(see \cite{Z96}, and also \cite{P00}, \cite{PS02},\cite{PS03}).
The simple formula (\ref{F0}) raises the intriguing possibility
that the behavior of the Monge-Amp\`ere operator $\O_{\Phi}^{n+1}$
and geodesics in ${\cal H}$ can be linked even more directly
to stability and constant scalar curvature K\"ahler metrics.

\bigskip

$\bullet$ Fix now two potentials $h_0$ and $h_1$ in
${\cal H}$,
and let $ {\underline s}^{(0)}$, $ {\underline s}^{(1)}$
be bases for $H^0(L^k)$ which are orthonormal
with respect to the metrics induced by $h_0$ and $h_1$.
As in the Introduction, let $\sigma\in GL(N_k+1)$ be the
corresponding change of bases, $\sigma\cdot {\underline s}^{(0)}
= {\underline s}^{(1)}$. We may assume that $\sigma$ is diagonal
with eigenvalues $e^{\lambda_0},\cdots,e^{\lambda_{N_k}}$.
We consider the functions $\phi(t;k)$ defined as in
(\ref{phi}).

\bigskip
\begin{lemma}
\label{simple bounds} Let $\O_0 = \pi_1^*\o_0 $
where $\o_0=R(h_0)$.

{\rm (a)} For each $k$, let $\Phi(k)$
be the extension of $\phi(t;k)$ to
$\bar M=X\times A$ as in (\ref{Phi}),
and let $\Omega_{\Phi(k)}$ is the corresponding $(1,1)$-form
as in (\ref{OPhi}). Then $\Omega_{\Phi(k)}$ is
a smooth positive $(1,1)$-form,
\be
\label{positive}
\Omega_{\Phi(k)}\geq 0.
\ee
In particular, the $(n+1,n+1)$-form $\Omega_{\Phi(k)}^{n+1}$
form is a positive smooth measure on $\bar M$;

\medskip

{\rm (b)} There is a constant $C>0$ which
does not depend on $k$ such that
\be
C^{-1}k\,\leq {\rm max}_{0\leq j\leq N}|\lambda_j|\,
\leq\,C\,k;
\ee
{\rm (c)} With the same constant $C>0$, we also have
\be
|\phi(t;k)|+
|\dot\phi(t;k)|\ \leq \ 4\,C.
\ee
\end{lemma}

\bigskip

Proof. Locally on $\bar M$, we have $e^{2\lambda_j
t}=|w^{\lambda_j}|^2$, where $w^{\lambda_j}=e^{\lambda_j(t+is)}\in A$
is a holomorphic function. Since the logarithms of sums of squares
of absolute values of holomorphic functions are plurisubharmonic in 
the
usual sense on ${\bf C}^{n+1}$,
and since we can write
\bea
\Omega_{\Phi(k)}
&=&
\Omega_0+{\sqrt{-1}\over 2}\ddb({1\over
k}\log\sum_{j=0}^{N_k}e^{2\lambda_j t}|s_j(z)|_{h_0^k}^2)
\nonumber\\
&=&
{\sqrt{-1}\over 2k}\ddb \log \sum_{j=0}^{N_k}|w^{\lambda_j}
  s_j(z)|^2,
\eea
where $|s_j(z)|$ is just the absolute value of $\hat s_j(z)$
in a local trivialization, the positivity
of $\Omega_{\Phi(k)}$ follows. Here we have simplified the
notation by denoting $s_j^{(0)}$ just by $s_j$.

\medskip
We turn next to the proof of (b). Let $\phi=\log {h_1\over h_0}$,
and order the eigenvalues $\lambda_j$ so that
$\lambda_0\geq\lambda_1\geq\cdots\geq\lambda_N$. Clearly,
\be
{2\over k}\lambda_N
\leq
{1\over k}
\log\,
{\sum_{j=0}^{N_k}e^{2\lambda_j}|s_j(z)|^2
\over
\sum_{j=0}^{N_k}|s_j(z)|^2}\leq {2\over k}\lambda_0.
\ee
By the Tian-Yau-Zelditch theorem, the expression in the
middle tends to $\phi$ as $k\to\infty$. Thus,
if we set $C_1={\rm sup}\,\phi$ and $C_2=-{\rm inf}\,\phi$,
we have $\lambda_0\geq C_1k/4$ and $\lambda_N\leq -C_2k/4$, for
$k$ large enough.

\medskip
To get inequalities in the opposite direction, we note that
for each $j$,
\be
\I |s_j|^2_{h_0^k}{\o_0^n\over
n!} = 1,
\qquad
\I {1\over N_k}\s_{l=0}^N |s_l|_{h_0^k}^2{\o_0^n\over
n!}=1.
\ee
This implies that for some $x\in X$,
we have $|s_j|^2 \geq {1\over N_k}\s_{l=0}^N |s_l|^2$. Choosing
$j=0$, we obtain for
$k$ large,
\be
{1\over k}\log {e^{\l_0}\over N_k }\
\leq \
{1\over k}\log{e^{\l_0 }|s_0(x)|^2\over |\us(x)|^2 }\ \ \leq \
{1\over k}\log{\s e^{\l_l}|s_l(x)|^2\over
   \s |s_l(x)|^2}\
\ \leq\ 2C_1.
\ee
This shows that ${\l_0}\leq {2C_1k}+\log N_k$. But $N_k\sim k^n$. Thus
we conclude
\be
{C_1\over 2}k \ \leq \l_0 \ \leq 3C_1 k,
\qquad k>>1.
\ee
To get the bound on $\l_{N_k}$, interchange the roles of $h_0$ and
$h_1$,
which changes $\phi$ to $-\phi$ and $\sigma$ to $\sigma^{-1}$. Thus we
have
\be
{C_2\over 2}k \ \leq \ -\l_{N_k} \ \leq\  3C_2 k,
\ee
and (b) is proved.

\medskip
The first inequality in (c), namely that $|\phi(t;k)|$ is uniformly
bounded, follows immediately from the bounds
in (b) for $|\lambda_j|$. To get the second inequality,
we write
\be
\label{phidot}
\dot\phi(t;k)
= {1\over k}{\sum_{j=0}^{N_k}2\lambda_je^{2\lambda_j t}|s_j(z)|^2\over
\sum_{j=0}^{N_k}e^{2\lambda_jt}|s_j(z)|^2}
\ee
and hence
\be
|\dot\phi(t;k)|\leq 2\,{{\rm max}|\lambda_j|\over k}.
\ee
Applying (b) again, we find that $|\dot \phi|$ is
uniformly bounded. Q.E.D.

\bigskip

$\bullet$ We can now state and prove the following theorem,
which is a key step in the construction of a generalized
solution of the Monge-Amp\`ere equation $\Omega_\Phi^{n+1}=0$:

\medskip

\begin{theorem}
\label{volume decay}
The volume of $X\times A$ approaches
zero as $k$ tends to infinity. More
precisely,
\be\label{volume decay two}
0\ \leq\ \I_{X\times A}\O_{\Phi(k)}^{n+1}\ \leq\
{C\over k}
\ee
where $C$ is a constant which is
independent of $k$. In particular, the positive measures
$\Omega_{\Phi(k)}^{n+1}$ tend weakly to 0.
\end{theorem}

\bigskip
{\it Proof.} Recall that $\dot\phi(t;k)$ has been evaluated and is
given by (\ref{phidot}). In particular, it is independent of the
choice of metric on $L^k$.
We now apply (\ref{volume}) to the path
of Fubini-Study metrics $\phi(t;k)$ defined
by (\ref{phi}) and obtain
\be
\label{volume one}
\I_{X\times A}\O_{\Phi(k)}^{n+1}\ = \
{2\over k}\I_X { \s_{j=0}^{N_k} \l_j| s_j^{(1)}|^2\over
\s_{j=0}^{N_k}| s_j^{(1)}|^2}\ \o_{\phi(1)}^n\ - {2\over k}\I_X
{ \s_{j=0}^{N_k}
\l_j| s_j^{(0)}|^2\over
\s_{j=0}^{N_k}| s_j^{(0)}|^2}\ \o_{\phi(0)}^n\
\ee
Here $\Phi(k)$ is defined by (\ref{Phi})
with $\phi(t)$ replaced by $\phi(t;k)$.
Since $h_1(k)=h_{\hat \us^{(1)}}$
the formula (\ref{fubini two}) implies
$$k^{-n}
\s_{j=0}^{N_k} | s_j^{(1)}|_{h_1^k(k)}^2
\ = \
\s_{j=0}^{N_k} |\hat s_j^{(1)}|_{h_1^k(k)}^2\ = \ 1
$$
and, observing that $\o_{\phi(1)}=\o_1(k)$ and
$\o_{\phi(0)}=\o_0(k)$,
 we can rewrite (\ref{volume one}) as
\be
\label{volume two}
\I_{X\times A}\O_{\Phi(k)}^{n+1}\ = \
{2\over k^{n+1}}\I_X { \s_{j=0}^{N_k} \l_j
| s_j^{(1)}|_{h_1^k(k)}^2}
\ \o_1^n(k)\ - {2\over k^{n+1}}\I_X
{ \s_{j=0}^{N_k}
\l_j| s_j^{(0)}|_{h_0^k(k)}^2}
\ \o_0^n(k).\
\ee
Now observe that
$$
\I_X |s_j^{(1)}|_{h^k_1(k)}^2\ \o^n_1(k)\ = \ \I_X
|s_j^{(1)}|_{h^k_1}^2\
\o_1^n\cdot {h^k_1(k)\over h^k_1}\cdot {\o_1^n(k)\over \o_1^n}.
$$
On the other hand, (\ref{powers of k}) implies ${h^k_1(k)\over
h^k_1}=1+O({1\over k})$ and ${\o^n_1(k)\over \o^n_1}=1+O({1\over 
k^2})$.
Moreover,
$\I_X |s_j^{(1)}|^2_{h_1^k}\ \o_1^n\ = \ 1 $
since the $s_j^{(1)}$ are orthonormal with
respect to $h_1^k$.
Thus
$$
\I_X |s_j^{(1)}|_{h^k_1(k)}^2\ \o^n_1(k)\ = \
1\ +\ O({1\over k})
$$
and hence
\be\label{big O}
{2\over k^{n+1}}\I_X { \s_{j=0}^{N_k} \l_j
| s_j^{(1)}|_{h_1^k(k)}^2}
\ \o_1^n(k)\\ - \
{2\over k^{n+1}}{ \s_{j=0}^{N_k} \l_j
}
\ = \
O({1\over k^{n+2}})\cdot N_k\cdot\,\max_{0\leq j\leq N_k}
|\l_j|
\ee

Now the Riemann-Roch theorem implies
$N_k=O(k^n)$ and thus we obtain from
Lemma \ref{simple bounds} that the
right side of (\ref{big O}) is of
the size $O({1\over k})$. Theorem \ref{volume decay}
follows now from (\ref{volume two}) and
(\ref{big O}). Q.E.D.

\section{Generalized solutions of the Monge-Amp\`ere equation}
\setcounter{equation}{0}

In the previous section, we have seen that the functions $\Phi(k)$ 
form
a uniformly bounded sequence of
functions whose Monge-Amp\`ere operators $\Omega_{\Phi(k)}^{n+1}$
are positive measures on $\bar M$ which
tend to~0. The Chern-Levine-Nirenberg inequality \cite{CLN}
implies that if any subsequence of the
$\Phi(k)$'s converges uniformly, then its limit $\Phi$
would satisfy
the Monge-Amp\`ere equation $\Omega_\Phi^{n+1}=0$ in the
generalized sense. A major problem is the fact
that the bounds available to us at the present time (c.f. Lemma 1) are
not strong enough to guarantee the existence of a uniformly convergent
subsequence of the $\Phi(k)$'s. Of course, weakly convergent
subsequences can always
be found. However, it is well-known that the Monge-Amp\`ere operator
$PSH\cap C^0\ni\Phi\to \Omega_\Phi^{n+1}$ is not lower semi-continuous
under weak limits \cite{L}.

\medskip
To circumvent these difficulties, we shall formulate and establish
extensions of the classical convergence and uniqueness
theorems of Bedford-Taylor for the Monge-Amp\`ere operator
for domains in ${\bf C}^n$ to the case of K\"ahler manifolds
with boundary.

\medskip

\subsection{Convergence of approximate solutions}

The first of these extensions is the following convergence theorem,
where the key hypothesis is the existence of uniform $C^0$ bounds
for some transversal derivative of the $\Phi(k)$'s at the boundary of
the manifold $\bar M$.

\medskip

Let $\bar M$ be a compact complex manifold with smooth
boundary and  $M\sub  \bar M$ be the
interior of $\bar M$. Let $\bar M=\cup_{\alpha=1}^NU_\al$
be a covering of $\bar M$ by a finite number of coordinate
charts $U_\al$.
Fix a smooth closed $(1,1)$-form $\O_0$ on $M$,
and let $\Psi_\al$ be smooth potentials
for $\O_0$ on $U_\al$, that is,
$\O_0={\sqrt{-1}\over 2}\ddb\Psi_\al$ on $U_\al$. We define the class 
of
$\O_0$-plurisubharmonic functions on $M$
to be the class of functions $\Phi$ on $\bar M$
by
\be
PSH(M,\O_0)=\{\Phi;\ \Psi_\al+\Phi\  is\ plurisubharmonic\ on\ U_\al,
\ 1\leq\al\leq N\}.
\ee
Note that this condition means that $\Psi_\al+\Phi$ is upper 
semi-continuous
and satisfies the sub-mean value property.
However, it is a stronger condition than the condition
$\Omega_0+{\sqrt{-1}\over 2}\ddb\Phi\geq 0$ by itself, since it is a 
pointwise condition,
while the condition $\Omega_0+{\sqrt{-1}\over 2}\ddb{\Phi}\geq 0$ 
depends only on the values
of $\Phi$ almost everywhere.

\bigskip

Next, let $\phi$ be a continuous function on
$\pl\bar M$, and $\O_0$ a real, smooth
closed $(1,1)$ form on $M$. Set $d=\pl+\bar\pl$,
$d^c={\sqrt{-1}\over 4}(\bar\pl-\pl)$,
so that $dd^c={\sqrt{-1}\over 2}\ddb$.

\bigskip

\begin{definition}
\label{dirichlet}
Let $\Phi:\bar M\ra\R$ be an
upper-semicontinuous function. We say $\Phi$
is a solution of the Dirichlet problem with
boundary values $\phi$ if

\begin{enumerate}

\item $\Phi\in PSH(M,\O_0) $
\item $\Phi$ is continuous at $p$ for all
$p\in\pl\bar M$ and $\Phi|_{\pl\bar M}=\phi$.
\item $(\O_0+dd^c\Phi)^{m}=0 $ on $\bar M$ where
$m=\dim(M)$.

\end{enumerate}
\end{definition}

\v

We wish to obtain a solution of the
Dirichlet problem on $\bar M $ from a sequence
of approximate solutions.
Let $Y$ be a smooth real nowhere vanishing
vector field
on a neighborhood $U\sub \bar M$ of
$\pl\bar M$
which is transversal to $\pl\bar M$,
in the sense that for $p\in\pl\bar M$,
the vector $Y(p)$ is not tangent to $\pl\bar M$.
\v

\begin{theorem}
\label{approximate} Assume
$\O_0^{m}=0
$. Let $\Phi_k\in PSH(M,\O_0)\cap C^\i(\bar M)$
have the properties:

\begin{enumerate}
\item $|\Phi_k|$ is uniformly bounded on $\bar
M$.
\item $|Y(\Phi_k)| $ is uniformly bounded on
$U$;
\item $\Phi_k|_{\pl\bar M} \ra \phi$ uniformly;
\item There is a sequence
$a_k > \sup_{x\in\pl\bar M}
|\Phi_k(x)-\phi(x)|$ with the property
$a_k\searrow 0$ (i.e. $a_k$ is decreasing and
approaching zero) and
$\sum_ka_k <\i$.
\item $\lim_{k\to\i}\I_{M}\ (\O_0+dd^c\Phi_k)^m  =
0$

\end{enumerate}

Let
$$  \Phi\ = \ \lim_{k\to\i}\big[\
\sup_{l\geq k} \Phi_l\ \big]^*
$$
Then $\Phi$ is a solution of the Dirichlet
problem with boundary values $\phi$.
\end{theorem}

\v

{\it Remark}: Assumption 4. implies Assumption 3.
On the other hand, after passing to a
subsequence, Assumption 3. implies Assumption 4.
\v

{\it Proof.}  We
want to choose a sequence $c_k\searrow 0$ with
the  property in such a way that $\Phi_k+c_k$
is monotonically decreasing on $\pl M$. To do
this, we just define $c_k=2\s_{j\geq k}a_j $.
Then
$c_{k}-c_{k+1}=2a_k$. Moreover, if $x\in\pl M$,
then
\be 
\label{wma} 
\Phi_k(x)-\Phi_{k+1}(x)\ >
(\phi-a_k)- (\phi+a_{k+1}) \ \geq -2a_k\ = \
c_{k+1}-c_k
\ee
Thus, replacing $\Phi_k$ by $\Phi_k+c_k$, we
may assume that
$\Phi_k|_{\pl\O}>\Phi_{k+1}|_{\pl\O}$.
\v
Let
$$ W_k\ = \ \big[\ \sup_{l\geq k} \Phi_l\
\big]^*
$$
Then $W_k\in PSH(M,\O_0)$ by Theorem 5.7 of
\cite{D}. Moreover, (\ref{wma}) implies
that $W_k=\Phi_k$ in an open neighborhood of
$\pl\bar M$. Thus
\be \label{boundary}
   \I_{M} (\O_0+dd^cW_k)^m\ = \  \I_{M}
(\O_0+dd^c\Phi_k)^m
\ee

This follows from the following simple lemma:
\v
\begin{lemma} \label{stokes0}
{\it\
Let $W,\Phi\in PSH(M,\O_0)$ with $\O_0^m=0$ and
assume that $W=\Phi $ on some neighborhood of $\pl\bar
M$. Then
$$  \I_{M} (\O_0+dd^cW)^m\ = \  \I_{M}
(\O_0+dd^c\Phi)^m
$$

}
\end{lemma}

{\it Proof of Lemma \ref{stokes0}}.
To see this, let
$K$ be any compact subset such that $W=\Phi$
on $M\backslash K$, and let $\Psi\in C_0^\infty(M)$
be a function with $\Psi=1$ on a neighborhood of $K$.
Then, expanding
$(\O_0+dd^cW)^m$, we obtain
$$ \I_M \Psi(\O_0+dd^cW)^m\ =
\I_M \Psi\O_0^m\ + \ \I_M\Psi dd^c(\Theta_W)\ = \
\I_{M\backslash K} dd^c\Psi\wedge \Theta_W
$$
Here $\Theta_W=\s_{k=0}^{m-1}
\big({m\atop k}\big) \O_0^k\wedge W\wedge
(dd^cW)^{m-1-k}
$.
Since $\Theta_W=\Theta_\Phi$ on
$M\backslash K$, the lemma is proved.

\v
Now assumption 5. of the theorem together
with (\ref{boundary}) implies that
$(\O_0+dd^cW_k)^m\ra 0 $ weakly. On the other
hand, since $W_k\ra \Phi$ monotonically, the
Bedford-Taylor monotonicity theorem (Theorem 2.1
of \cite{BT82}; this is stated for domains in
$\C^n$, but generalizes in a straightforward
fashion to manifolds) we have
$(\O_0+dd^cW_k)^m\ra (\O_0+dd^c\Phi)^m
$ weakly. Thus we have  $(\O_0+dd^c\Phi)^m=0$.
To finish the
proof of the theorem, we must show that $\Phi$
is continuous at the boundary and has the right
boundary values.
\v
Let $\e>0$. Choose $k_0$
such that $k\geq k_0 \Longrightarrow
\sup_{x\in\pl\bar M}|\Phi_k(x)-\phi(x)|<\e$.
Extend
$\phi$ to a continuous function on a
neighborhood
$U\sub \bar M$ of $\pl\bar M$ in such a way
that
$\phi$ is constant on the flow lines of $Y$.
Then assumption 2. implies that if $U$
is sufficiently small, then
\be
\sup_{x\in U}|\Phi_k(x)-\phi(x)|<2\e
\ee
and thus $\sup_{x\in U}|\Phi(x)-\phi(x)|\leq 2\e$.
In particular, $\Phi|_{\pl\O}=\phi$  and
$\Phi$ is continuous at all points
$p\in\pl\bar M$. This proves the theorem.
Q.E.D.

\subsection{The domination principle for the Monge-Amp\`ere operator}

The proof of the extension of the Bedford-Taylor uniqueness
theorem in this section follows
closely the original arguments of
\cite{BT76, BT82}, and especially the exposition of Blocki
\cite{B}. 

\v

When we consider generalized solutions of a
partial differential equation, a particularly desirable
property is their uniqueness. For bounded domains in ${\bf
C}^n$, the uniqueness of the generalized solution of the
Dirichlet problem for the Monge-Amp\`ere equation
in the class $PSH\cap L^{\infty}$ has been established by
Bedford and Taylor \cite{BT76}. It seems that this uniqueness
theorem should extend as well to bounded domains in
K\"ahler manifolds, at least if good smooth approximations of
$\Omega_{\Phi(0)}$-plurisubharmonic functions exist.
Although there are now many powerful approximation theorems
(see \cite{D, GZ} and references therein), we found it more
convenient to extend the Bedford-Taylor uniqueness theorem
to a situation adapted to the problem at hand, in the
spirit of the earlier extension. The key hypothesis which
we will exploit is a capacity zero condition.

\bigskip
Recall the following notion of capacity of a set introduced by
Bedford and Taylor \cite{BT82}: If $E\sub U$ is a Borel
subset of a bounded domain $U\sub \C^n$ then
\be
\label{capacity}
c(E,U)
=
{\rm sup}\bigg\{\int_E(dd^c v)^n; v\in PSH (U),
0\leq v \leq 1\}.
\ee
For our purposes, we can adapt this notion to K\"ahler manifolds
as follows: Let $M=\cup_{\al=1}^N U_\al$ be a finite cover of $M$
by coordinate neighborhoods. Then we say $c(E,M)<\e$ if
we can write $E=\cup_\al E_\al$ with $E_\al\sub U_\al$ a Borel
subset and 
\be
\label{capacityepsilon}
\s_\al c(E_\al, U_\al)\ < \ \e.
\ee
We say that $c(E,M)=0$ if $c(E,M)<\epsilon$ for every
$\epsilon>0$.

\begin{lemma}\label{trivial}
There is a constant $C>0$ with the
following property: If
$E\sub M$ is a Borel subset, $\e>0$  and
$\Phi\in PSH(\O_0,M)$, then
\be
\I_E (\O_0+dd^c\Phi)^m\ \leq \ C\e(1+\sup |\Phi|)^m
\ee
if $c(E,M)<\e$. In particular, if $c(E,M)=0$, then
for all functions $\Phi\in PSH(\O_0,M)$, we have 
\be
\I_E(\O_0+dd^c\Phi)^m=0.
\ee
\end{lemma}

{\it Proof.} Fix a smooth potential $\Psi_\al$ on $U_\al$ such
that $\O_0=dd^c\Psi_\al$. Then
\bea
\I_E (\O_0+dd^c\Phi)^m\ &\leq& \
\s_\al \I_{E_\al} (\O_0+dd^c\Phi)^m\ = \
\s_\al \I_{E_\al} (dd^c(\Psi_\al+\Phi))^m\ 
\nonumber\\
&\leq& \
\s_\al(\sup|\Psi_\al+\Phi|)^m\,c(E_\al,U_\al).
\eea
Since the $\Psi_\al$ are fixed, we have
$|\Psi_\al+\Phi|\leq C(1+|\Phi|) $ and the lemma follows.

\bigskip
The following lemma follows immediately from the
quasi-continuity theorem of Bedford-Taylor
\cite{ BT82}:

\bigskip

\begin{lemma}\label{qc}
Let $\Phi\in PSH(M,\O_0)$. Then for every
$\e>0$, there is an open
set
$G\sub M$ such that $c(G,M)<\e$ and $\Phi$ is continuous
on $M\backslash  G$. 
\end{lemma}
\v
{\it Proof.}
Since $\Psi_\al+\Phi$ is a plurisubharmonic 
function on $U_\al$, the quasi-continuity theorem of
Bedford-Taylor (\cite{BT82}, Theorem 3.5)
implies that there is an open set $G_\al\sub U_\al$
such that $\Psi_\al+\Phi$ is continuous on $U_\al\backslash G_\al$
and $c(G_\al,U_\al)<\e$. Let $G=\cup_\al G_\al$. Then $\Phi$
is continuous on $M\backslash G$ and, by definition,
$c(G,M)<N\,\e$.
\v

\v
We also require a notion of ``nearly continuous" functions:

\bigskip
\begin{definition}
\label{nearlycontinuous}
We say that a bounded function
$v:\bar M\to\R$ is ``nearly continuous" if
\begin{enumerate}
\item There exists a lower semi-continuous
function $v_0$ on $\bar M$ such that
$v=v_0^*$;
\item $\{ v_0<v\} $ has capacity zero,
that is, $c(\{v_0<v\},M)=0$;
\item $v=v_0$ on $\pl M$.
\end{enumerate}

\end{definition}

\bigskip

With this notion, we shall prove the following:

\bigskip

\begin{theorem}
\label{dom}
Assume that $\bar M$ is a complex manifold of dimension
$m$ with smooth boundary, let $M=\bar M\backslash \pl M$
and let
$\O_0$ be a real closed smooth $(1,1)$ form on
$M$ satisfying $ \O_0^m=0$.  Let
$u,v\in PSH(M,\O_0)\cap L^\i$ be such that $(u-v)_*\geq
0$ on
$\pl \bar M$. Assume as well:

\begin{enumerate}

\item $u$ is continuous;
\item There is a decreasing sequence $v_k$ of nearly
continuous functions in $PSH(M,\O_0)$ such that  $v_k\searrow
v$;
\item For every $\d>0$ there is a compact set $K\sub M$
such that $v_k|_{M\backslash K} < v|_{M\backslash
K}+\d$ for $k$ sufficiently large.

\end{enumerate}

Then
\be\label{domination} \I_{u<v}(\O_0+dd^cv)^m\ \leq
\I_{u<v}(\O_0+dd^cu)^m.
\ee
\end{theorem}

\v
{\it Proof.} We divide the proof into several steps.

\medskip

Step 1. If we replace $u$ by $u+\d$, then we have
$\{u+\d<v\}\uparrow \{u<v\}  $ as $\d\downarrow
0$. Since for any positive measure $\m$ we
have $\m(E_j)\ra \m(\cup E_j)$ whenever $E_j$ is
an increasing family of measurable sets, we may
replace $u$ by $u+\d$. Thus we may assume that
\be
(u-v)_*\geq \d \hbox{ and that $M'=\{u<v\} $\ 
is relatively compact in $M$}.
\ee
\v
Step 2. We prove the theorem under the assumption that
$v$ is continuous (in which case the hypotheses 2. and 3. are
automatic, since we can take $v_k=v$ for all $k$).
\v
For $\e>0$ let $u_\e=\max(u+\e,v)$. Then  $u_\e= u+\e$
on a neighborhood of $\pl M'$. Since $dd^c
(u+\e)= dd^c u$, we can invoke Lemma 2 and conclude that
\be
\I_{M'}(\O_0+dd^cu_\e)^m\ =
\I_{M'}(\O_0+dd^cu)^m.
\ee
On the other hand,
$u_\e\downarrow
v$ on $M'$ and, by the Bedford-Taylor
monotonicity theorem,
$$ (\O_0+dd^cu_\e)^m \ \ra\ (\O_0+dd^cv)^m
\ on\ M'
\ \ \ \hbox{(weak convergence of measure)}.
$$
Since $\{u<v\}$ is open, we obtain
$$ \I_{\{u<v\}} (\O_0+dd^cv)^m\ \leq \
\liminf_{\e\to 0} \I_{\{u<v\}} (\O_0+dd^c
u_\e)^m.
$$
This completes step 2.

\v
Step 3. Now we treat the  case where $v$ itself is nearly
continuous (in which case the hypotheses 2. and 3. are
again automatic, since we can take $v_k=v$ for all $k$). 
This step will be parallel to the argument in Step 2:
\v
For $\e>0$ let $u_\e=\max(u+\e,v)$. Then $u_\e\searrow v$
on the open set $\{u<v_0\}$ (this set is open since
$v_0$ is lower semi-continuous). In particular, the
Bedford-Taylor monotonicity theorem implies
$(\O_0+dd^cu_\e)^m\ra (\O_0+dd^cv)^m$ weakly 
on $\{u<v_0\}$ (as measures
or as currents - the two notions of weak convergence
are equivalent).
Thus we have
$$ \I_{\{u<v\}}(\O_0+dd^cv)^m\ = \
\I_{\{u<v_0\}}(\O_0+dd^cv)^m
\ \leq \ \liminf_{\e\to 0}\I_{\{u<v_0\}}(\O_0+dd^c u_\e)^m
$$
$$
\ = \liminf_{\e\to 0}\I_{\{u<v\}}(\O_0+dd^c u_\e)^m
$$
The first equality follows from the  assumption that
$\{v_0<v\} $ has capacity zero and the inequality from
the fact that $(\O_0+dd^cu_\e)^m\ra (\O_0+dd^cv)^m$ weakly. Here
we are making strong use of the fact that $\{u<v_0\} $
is open.
\v
Next we claim that for all $\e>0$:
\be\label{epsilon independence}
\I_{\{u<v\}}(\O_0+dd^c u_\e)^m\
=
\
\I_{\{u<v\}}(\O_0+dd^c u)^m.
\ee
To see this, let $A$ be an open set containing
$\{u<v\} $, which is relatively compact in $M$. Then, $p\in \pl A$ 
implies that
$p\notin A$, and we have
$u(p)\geq v(p)$, so $u_\e(p)=u(p)+\e>v(p)$. Since
the set $\{u+\e>v\}$ is open ($u $ is continuous and $v$
is upper semi-continuous), we see that
$u+\e=u_\e$ in a neighborhood of $\pl A$. Thus Lemma
\ref{stokes0} implies

$$
\I_A(\O_0+dd^c u_\e)^m\
=
\
\I_A(\O_0+dd^c (u+\e))^m\ = \ \I_A(\O_0+dd^c u)^m
$$
Since $(\O_0+dd^c u_\e)^m$ and $(\O_0+dd^c u)^m$ are positive
Borel measures which are finite on compact subsets
of $M$, they are both regular.
In particular
\bea
\I_{\{u<v\}}(\O_0+dd^c u_\e)^m\ &=& \ \inf_{A}
\I_{A}(\O_0+dd^c u_\e)^m
\nonumber\\
\I_{\{u<v\}}(\O_0+dd^c u)^m\ &=& \ \inf_{A}
\I_{A}(\O_0+dd^c u)^m
\eea
where the inf is taken over all open sets $A$ which
contain $\{u<v\} $.
This proves (\ref{epsilon independence}).

\v

Step 4. We treat the general case. From
step 1. we can assume that
$(u-v)_*\geq
\d>0$. Thus
$u>v_k$ on a  neighborhood of $\pl M$ for $k$
sufficiently large. By Lemma \ref{qc}, there is an open set
$G\sub M$ such that
$c(G,M)<\e$ (the capacity of
$G$ in $M$) and such that
$v$ is continuous on $F=M\backslash G$. Choose $\phi$
which is continuous $M$ such that $\phi=v$ on $F$. Then
$\{u<v\}\sub \{u<\phi\}\cup G $ implies
\bea
\I_{\{u<v\}}(\O_0+dd^cv)^m\ &\leq& \
\I_{\{u<\phi\}}(\O_0+dd^cv)^m\ + \
\I_G(\O_0+dd^cv)^m\nonumber\\
 &\leq& \
\I_{\{u<\phi\}}(\O_0+dd^cv)^m+C\e
\eea

where the last inequality follows from Lemma \ref{trivial},
the fact that $c(G,M)<\e$ and that the function $v$ is
bounded. Here the constant $C$ depends only on
the sup norm of $v$.
\v
Since $\{u<\phi\}$ is open and $v_k\searrow v$, 
the Bedford-Taylor monotonicity
theorem implies
$$ \I_{\{u<\phi\}}(\O_0+dd^cv)^m\ \leq \
\liminf_{k\to\i}
\I_{\{u<\phi\}}(\O_0+dd^cv_k)^m\
$$

Now
$$ \{u<\phi\}\ \sub \ \{u<v\}\cup G \ \sub\
\{u<v_k\}\cup G
$$
Thus

$$\I_{\{u<v\}}(\O_0+dd^cv)^m\ \leq\
\liminf_{k\to\i}
\I_{\{u<v_k\}}(\O_0+dd^cv_k)^m\ \ + \ 2C\e,
$$
since our assumptions imply that the functions $v_k$ are uniformly
bounded.
Next, the assumption 3. implies that the sets $\{u<v_k\}$ are all
contained in a relatively compact subset of $M$. 
Using Step 3,
$$\I_{\{u<v_k\}}(\O_0+dd^cv_k)^m\ \leq \
\I_{\{u<v_k\}}(\O_0+dd^cu)^m\
$$
and since $v_k\searrow v $, $\cap_k\{u<v_k\}=\{u\leq v\}$,
we can conclude that
$$\I_{\{u<v\}}(\O_0+dd^cv)^m\ \leq\
\I_{\{u\leq v\}}(\O_0+dd^cu)^m\ \ + \ 2C\e
$$

Since $\e$ is arbitrary:
$$\I_{\{u<v\}}(\O_0+dd^cv)^m\ \leq\
\I_{\{u\leq v\}}(\O_0+dd^cu)^m\
$$
Applying this last inequality to $u+\eta$ and $v$, for
some $\eta>0$:
$$\I_{\{u+\eta<v\}}(\O_0+dd^cv)^m\ \leq\
\I_{\{u+\eta\leq v\}}(\O_0+dd^cu)^m\
$$
Finally, taking the limit as $\eta\ra 0$ and noting
that $\cup_{\eta>0}\{u+\eta<v\}
=\cup_{\eta>0}\{u+\eta\leq v\}=\{u<v\}$, we obtain
(\ref{domination}).

\v

Next we prove another version of the domination
theorem, but this time with the roles of $u,v$ reversed:

\begin{theorem}
\label{domination two}
Assume that $\bar M$ is a complex manifold of dimension
$m$ with smooth boundary, let $M=\bar M\backslash \pl M$
and let
$\O_0$ be a real closed smooth $(1,1)$ form on
$M$ with the property: $ \O_0^m=0$.  Let
$u,v\in PSH(\bar M,\Omega_0)\cap L^\i$ satisfy $(u-v)_*\geq
0$ on
$\pl \bar M$. Assume as well:

\begin{enumerate}

\item $v$ is continuous;
\item There is a decreasing sequence $u_k$ of nearly
continuous PSH functions on $M$ such that  $u_k\searrow
u$.

\end{enumerate}

Then
\be\label{domination} \I_{u<v}(\O_0+dd^cv)^m\ \leq
\I_{u<v}(\O_0+dd^cu)^m.
\ee

\end{theorem}

{\it Proof}. The proof is parallel to
that of Theorem \ref{dom}, although there are some important
differences. We again divide it into several steps.

\v
Step 1.
As before, we may assume that
$$(u-v)_*\geq \d \hbox{ and that $M'=\{u<v\}$
is relatively compact in $M$}
$$
Moreover, we have
$$\hbox{  $\{u_0<v\}$
is relatively compact in $M$}
$$
To see this, observe first that $u\geq v+{\d/2}$
in an open neighborhood of $\pl M$. Also,
$\{u_0-u>-\d/4\}$ is open since $u_0$ is lower
semi-continuous and $u$ is upper semi-continuous.
Since $u_0=u$ on $\pl M$, the set 
$\{u_0-u>-\d/4\}$ is an open neighborhood of
$\pl M$. Thus $u_0>v+\d/4$ in an open neighborhood
of $\pl M$.
\v

Step 2. Now we treat the  case where $u$ itself
is nearly continuous (in which case hypotheses 2.
and 3. are automatic, since we can take $u_k=u$
for all $k$).
\v
For $\e>0$ let $u_\e=\max(u+\e,v)$. Then $u_\e\searrow v$
on the open set $\{u<v\}$ (this set is open since
$u$ is upper semi-continuous). As before, the
Bedford-Taylor monotonicity theorem implies
$(\O_0+dd^cu_\e)^m\ra (\O_0+dd^cv)^m$ weakly
and, making use of the fact that $\{u<v\}$
is open, we get
\bea
\I_{\{u<v\}}(\O_0+dd^cv)^m\
\ &\leq& \ \liminf_{\e\to 0}\I_{\{u<v\}}(\O_0+dd^c u_\e)^m
\nonumber\\
&=& \ \liminf_{\e\to 0}\I_{\{ u_0<v\}}(\O_0+dd^c
u_\e)^m,
\eea
where $u_0$ is the function which appears in the definition
of the near continuity of $u$.
Next we claim that for all $\e>0$:
\be
\label{epsilon independence2}
\I_{\{u_0<v\}}(\O_0+dd^c u_\e)^m\
=
\
\I_{\{u_0<v\}}(\O_0+dd^c u)^m
\ee
The proof is similar to that of (\ref{epsilon independence}),
using $A$ an open set containing
$\{u_0<v\} $, $u_0(p)\geq v(p)$ for  $p\in \pl A$
so that $u_0(p)+\e>v(p)$. We use now
the continuity of $v$ and the lower semi-continuity
of $u_0$ to deduce that
the set $\{u_0+\e>v\}$ is open, and
$u+\e=u_\e$ in a neighborhood of $\pl A$. As before, Lemma
2 implies
\be
\I_A(\O_0+dd^c u_\e)^m\
=
\
\I_A(\O_0+dd^c (u+\e))^m\ = \ \I_A(\O_0+dd^c u)^m,
\ee
and, using the fact that $(\O_0+dd^c u_\e)^m$ and $(\O_0+dd^c u)^m$ are 
positive
Borel measures which are finite on compact subsets
of $M$,
\bea
\I_{\{u_0<v\}}(\O_0+dd^c u_\e)^m\ &=& \ \inf_{A}
\I_{A}(\O_0+dd^c u_\e)^m\ 
\nonumber\\
\I_{\{u_0<v\}}(\O_0+dd^c u)^m\ &=& \ \inf_{A}
\I_{A}(\O_0+dd^c u)^m
\eea
where the inf is taken over all open sets $A$ which
contain $\{u_0<v\} $. This proves (\ref{epsilon independence2}).
Since $u$ and $u_0$ differ only on a set of capacity 0,
we obtain
the desired inequality.

\v

Step 4. We treat the general case: From
step 1. we can assume that
$(u-v)_*\geq
\d>0$. Thus
$u_j>v$ on a  neighborhood of $\pl M$ for $j$
sufficiently large. Choose an open set $G\sub M$ such
that
$c(G,M)<\e$ (the capacity of
$G$ in $M$) and such that
$u$ is continuous on $F=M\backslash G$. Choose $\phi$
which is continuous $M$ such that $\phi=u$ on $F$. Now
\bea
 \I_{\{u<v\}}(\O_0+dd^cv)^m\ &=& \
\lim_{j\to\i}\I_{\{u_j<v\}}(\O_0+dd^cv)^m\ \leq \
\lim_{j\to\i}\I_{\{u_j<v\}}(\O_0+dd^cu_j)^m
\nonumber\\
&\leq& \ \liminf_{j\to\i}\I_{\{u<v\}}(\O_0+dd^cu_j)^m
\eea
where we have made use of step 3. and the fact
that $u_j>v$ on $\pl M$ to prove
the first inequality.
Now we have
$\{u<v\}\sub \{\phi<v\}\cup G \sub (K\cap\{\phi\leq v\}\cup G$
where $K\sub M$ is a compact set such that $\{u<v\}\sub K$ so
\bea
\I_{\{u<v\}}(\O_0+dd^cu_j)^m\ &\leq& \
\I_{\{\phi\leq v\}\cap K}(\O_0+dd^cu_j)^m\ + \
\I_G(\O_0+dd^cu_j)^m\nonumber\\
&\leq& \
\I_{\{\phi\leq v\}\cap K}(\O_0+dd^cu_j)^m+C\e
\eea

where the last inequality follows from the fact
that $c(G,M)<\e$ and that the function $v$ is
bounded. Here the constant $C$ depends only on
the sup norm of $v$.
\v
Since $\{\phi\leq v\}\cap K$ is compact, the Bedford-Taylor
monotonicity theorem implies
$$ \limsup_{k\to\i} \I_{\{\phi\leq v\}\cap K}(\O_0+dd^cu_j)^m\
\leq \
\I_{\{\phi\leq v\}\cap K}(\O_0+dd^cu)^m\
$$

Now
$$
\{\phi\leq v\}\cap K \ \sub\ 
 \{\phi\leq v\}\ \sub \ \{u\leq v\}\cup G \
$$
Thus

$$\I_{\{u<v\}}(\O_0+dd^cv)^m\ \leq\
\I_{\{u\leq v\}}(\O_0+dd^cu)^m\ \ + \ 2C\e
$$

Since $\e$ is arbitrary:
$$\I_{\{u<v\}}(\O_0+dd^cv)^m\ \leq\
\I_{\{u\leq v\}}(\O_0+dd^cu)^m\
$$
Applying this last inequality to $u+\eta$ and $v$, for
some $\eta>0$:
$$\I_{\{u+\eta<v\}}(\O_0+dd^cv)^m\ \leq\
\I_{\{u+\eta\leq v\}}(\O_0+dd^cu)^m\
$$
Finally, taking the limit as $\eta\ra 0$, we obtain
(\ref{domination}).

\v

We can now state and prove the uniqueness theorem which
we need for the proof of Theorem 1:

\v

\begin{theorem}
\label{uniqueness}
Let $u,v$ be as in Theorem \ref{dom}, that is, 
\begin{enumerate}

\item $u$ is continuous;
\item There is a decreasing sequence $v_k$ of nearly
continuous PSH functions on $M$ such that  $v_k\searrow
v$;
\item For every $\d>0$ there is a compact set $K\sub M$
such that $v_k|_{M\backslash K} < v|_{M\backslash K}+\d$
for $k$ sufficiently large.

\end{enumerate}

Assume that
$(u-v)_*=(u-v)^*=0$ on $\pl M$ and that
$(\O_0+dd^cu)^m=(\O_0+dd^cv)^m=0$. Then
$u=v$.

\end{theorem}

{\it Proof.} We may assume, after replacing $u$ and $v$ by 
$u+C$ and
$v+C$ for some constant $C$, that $u$ and $v$ are
positive.
We wish to show $u=v$. Assume not: Let $\psi\in C^\i(M)$
be such that $\O_0+dd^c\psi>0$. Replacing $\psi$ by
$\psi-C$, we may assume that $\psi<0$ on $M$.
\v
Case 1.
$\{u<v\}\not=\emptyset$. This implies
$\{u<(1-\e)v+\e\psi\}\not=\emptyset
$ for some $\e>0$. Let $p\in \{u<(1-\e)v+\e\psi\}$ and
let $D$ be a disk in some coordinate neighborhood of $p$.
Then $D\cap \{u<(1-\e)v+\e\psi\}$ has non-zero Lebesgue
measure  (in general, if $u,v$ are psh functions such
that $u=v$ almost everywhere in a disk $D$, then
$u=v$ everywhere in $D$; this follows from local
regularization). Now we have, using
Theorem \ref{dom}
\bea
0\ = \  \I_{\{u<(1-\e)v+\e\psi\}}(\O_0+dd^cu)^m\ 
&\geq&
\I_{\{u<(1-\e)v+\e\psi\}}(\O_0+ dd^c[ (1-\e)v+\e\psi])^m
\nonumber\\
&\geq&\  \e^m\I_{\{u<(1-\e)v+\e\psi\}}
(\O_0+dd^c\psi)^m\ > 0
\eea

Case 2. $\{v<u\}\not=\emptyset$. This is treated
exactly in the same was as in case 1 except that
we use Theorem \ref{domination two} instead of
Theorem \ref{dom}. Q.E.D.

\section{Proof of Theorem 1}
\setcounter{equation}{0}

We can give now the proof of Theorem 1.

\medskip
First, we apply Theorem 3
to construct a generalized solution of the Dirichlet problem for the 
Monge-Amp\`ere equation
\be
\label{dboundary}
\Omega_\Phi^{n+1}=0\ on\ M=X\times A,
\qquad
\Phi|_{\partial M}=\phi,
\ee
where $\phi:\pl M\to {\bf R}$ is defined by $\phi|_{|w|=1}=0$ and
$\phi|_{|w|=e}=\log{h_0\over h_1}$. Define $\Phi(k)(z,w)=\phi(t;k)$,
where $t=\log |w|$ and $\phi(t;k)$ is defined as in
(\ref{phi}). Let $Y=\pl_t$. Then $|\Phi(k)|\leq C$
and $|Y(\Phi(k))|\leq C$ by Lemma 1.
We also have $||\Phi(k)|_{\pl\bar M}-\phi||_{L^\infty}\leq C{1\over 
k^2}$
by (\ref{powers of k}),
so that $\Phi(k)_{\pl\bar M}\to\phi$ uniformly.
Furthermore, Theorem 2 implies that $\int_M\O_{\Phi(k)}^{n+1}\to 0$
as $k\to\infty$. Thus, since $\sum_{k=1}^\i{1\over k^2}<\i$, we can 
apply Theorem 3,
and conclude that
\be
\Phi={\rm lim}_{k\to\infty}[{\rm sup}_{l\geq k}\,\Phi(l)]^*
\ee
is a generalized solution of the desired Dirichlet problem 
(\ref{dboundary}).

\medskip

Consider next the $C^{1,1}$ geodesic $\phi_t$ joining $\phi_0$ to
$\phi_1$ in the space ${\cal H}$ of K\"ahler potentials.
Let $\tilde\Phi(z,w)=\phi_t(z)$ with $t=\log|w|$
as before. We shall show that $\Phi=\tilde\Phi$.
To do this, we would like to apply Theorem 6 with
\be
u=\tilde\Phi,
\
(v_k)_0={\rm sup}_{l\geq k}\Phi(l),
\
v_k=(v_k)_0^*,
\
v=\Phi = \lim_{k\to\i} v_k
\ee
First, we show that $v_k$ is nearly continuous. Since $\Phi(l)$ is
smooth, $(v_k)_0$ is lower semi-continuous. Moreover,
$\{(v_k)_0<v_k\}$ has capacity 0, by Proposition 5.1 of
\cite{BT82}.
It remains to show that $v_k=(v_k)_0$ on $\pl\bar M$.
The equation (\ref{powers of k}) implies $|D\Phi(l)(z,w)|\leq C$,
for some constant $C$ independent of $l$, if $(z,w)\in\pl\bar M$
and $D$ is any derivative tangent to $\pl\bar M$.
Thus, 
if $\delta>0$,
there exists $\epsilon>0$ such that if $(z_0,w_0)\in\pl\bar M$ then
\be
\Phi(l)(z_1,w_0)+\d>\Phi(l)(z_0,w_0)>
\Phi(l)(z_1,w_0)-\delta,
\ee
for all $l$ if $|z_1-z_0|<\epsilon$. Also, $|Y(\Phi(l))|\leq C$
implies that 
\be
\Phi(l)(z_1,w_1)+\d>
\Phi(l)(z_1,w_0)>\Phi(l)(z_1,w_1)-\d\ {\rm if}\ |w_0-w_1|<\epsilon.
\ee
We have then 
$\Phi(l)(z_1,w_1)+2\delta>\Phi(l)(z_0,w_0)>\Phi(l)(z_1,w_1)-2\delta$,
which implies that 
\be
(v_k)_0(z_1,w_1)+2\d>
(v_k)_0(z_0,w_0)>(v_k)_0(z_1,w_1)-2\delta,
\ee
and hence
\be
\label{L1}
(v_k)_0(z_1,w_1)+2\d>
(v_k)_0(z_0,w_0)>v_k(z_1,w_1)-2\delta,
\ee
if $|z_1-z_0|<\e$, $|w_1-w_0|<\e$.
In particular,
$(v_k)_0(z_0,w_0)>v_k(z_0,w_0)-2\delta$ 
for all $\delta$ so
$(v_k)_0=v_k$ on $\pl M$ so  $v_k$ is indeed nearly continuous.

\medskip
Thus the first two assumptions of Theorem \ref{uniqueness}
are satisfied.

\medskip

In the next step, we will  need the bound
\be
\label{L2}
v(z_1,w_1)+2\d>
v(z_0,w_0)>v(z_1,w_1)-2\delta,
\ee
which follows by taking the limit of (\ref{L1}) as $k\to\i$.

\medskip

Next, we verify assumption 3: Since $v_k=(v_k)_0$ on
$\pl M$ we see that $v_k$ is both upper and lower
semi-continuous on $\pl M$ and thus $v_k$ is continuous
on $\pl M$. Moreover, by Theorem 3, $v=\phi$ on $\pl M$
so $v$ is also continuous on $\pl M$. Since $v_k\searrow v$
we see, by Dini's theorem, that $v_k\searrow v$ uniformly
on $\pl M$. Thus, for every $\d>0$ (\ref{L1}) and (\ref{L2}) imply 
that for $k>>0$

$$ v_k(z_1,w_1)-2\d <v_k(z_0,w_0)<v(z_0,w_0)+\d<v(z_1,w_1)+3\d
$$

for all $|z_1-z_0|<\e$ and $|w_1-w_0|<\e$.

\medskip

All the conditions of
Theorem~\ref{uniqueness}
are satisfied. We can thus conclude that
$\Phi=\tilde\Phi$.

\medskip
Finally, the uniform convergence of the functions
$[{\rm sup}_{k\geq l}\,\phi(t,k)]^*$ follows from
their upper semi-continuity and the compactness of $X$.
This is essentially Dini's theorem, and can be proven
as follows. Assume that $u_n$ is a sequence of upper
semi-continuous functions, decreasing to a
continuous limit $u$.
For each $\epsilon>0$,
the sets $\{x\in X;u_n(x)-u(x)<\epsilon\}$
form an open covering of $X$. Since $X$ is compact,
it admits a finite subcover, and since the sets are increasing
as $n$ increases, we must have $X=\cap_{n\geq N_{\epsilon}}
\{x\in X;u_n(x)-u(x)<\epsilon\}$ for some $N_\epsilon$.
Q.E.D.

\section{Remarks}
\setcounter{equation}{0}

We conclude with a few remarks.

\bigskip
$\bullet$ In \cite{D99}, Donaldson asks 
when two \K metrics
can be connected by a smooth geodesic. One way to approach
this problem is to establish a priori bounds on
the derivatives of $\phi(t;k)$. This was done  in
Lemma \ref{simple bounds}
  for $\phi$ and $\dot \phi$. Let
us now consider
$\ddot \phi$:
\be
\ddot\phi(0)\ = \ {1\over k}\s
(\l_\al-\l_\b)^2|s_i|_{h^k_0(k)}^2|s_j|_{h^k_0(k)}^2
\ee
Note that Lemma \ref{simple bounds} implies
$|\ddot\phi(0)|\leq Ck$, but this is not strong enough, since
we need a bound which is independent of $k$.

\v

Define a random variable $Z$ whose probability distribution
is given by
$$ P(Z=\l_\al)\ = \ |s_\al(z)|^2_{h^k_0(k)}
$$
This is indeed a distribution since the total
probability is $\s_\al |s_\al(z)|^2_{h^k_0(k)}=1$.
Moreover, $\s
(\l_\al-\l_\b)^2|s_i|_{h^k_0(k)}^2|s_j|_{h^k_0(k)}^2$ is
just  the variance of $Z$. Thus we need to prove that
the variance of $Z$ is bounded by $k$. In the simplest
case where $X=\P^1$ and the line bundle $L=O(1)$ and
the metric $h$ is the Fubini-Study metric, then an
easy computation shows that $Z$ is just the binomial
distribution with $k$ trials, where the probability $p$ of
flipping heads is a  function of $z\in \P^1$. As is
well known, the variance of the binomial distribution
is $kp$, which is the bound we need.
\v
In the case where $\o$ is the Fubini-Study metric
on $\P^1$, the eigenvalues of the change of basis
matrix are just $0,1,..., k$.
More generally, we can show that if $\phi$ is a
radially symmetric \K potential on $\P^1$, and
if $\l_0\leq \l_1\leq \cdots \leq\l_k$ are the
eigenvalues of the change of basis matrix, then
$|\l_j-\l_{j+1}| \leq C$ for some constant $C$,
independent of $k$. From this one can show without
difficulty that $\ddot\phi(0)$ is uniformly bounded.

\bigskip
$\bullet$ The function ${\rm lim}_{l\to\infty}[{\rm sup}_{k\geq
l}\,\phi(t;k)]^*(x)$ is equal almost everywhere to
the ${\rm lim\,sup}\phi(t;k)$. The convergence can also be
guaranteed to take place in a Sobolev norm of positive order.
Indeed, quite generally, the $L^2$-norm of $\pl \phi$
can be bounded by $||\phi||_{C^0}$
if $\phi$ is $PSH(X,\o_0)$-plurisubharmonic. Indeed,
\be
||\pl\phi||_{\o_0}^2
=\int_X \pl\phi\wedge\bar\pl\phi
\wedge\o_0^{n-1}
=
\int_X\phi\,\ddb\phi\wedge\o_0^{n-1}
=
\int_X\phi\,\o_\phi\wedge\o_0^{n-1}
-
\int_X\phi\,\o_0^n
\ee
and the right hand side can be bounded in turn by
\be
|\int_X\phi\,\o_\phi\wedge\o_0^{n-1}
-
\int_X\phi\,\o_0^n|
\leq
||\phi||_{C^0}(\int_X\o_\phi\wedge\o_0^{n-1}
+
\int_X\o_0^n)
=
2||\phi||_{C^0}\int_X\o_0^n.
\ee
(In fact, the same argument gives the following useful inequality
\be
J_{\o_0}(\phi)\leq 2n\,||\phi||_{C^0},
\ee
where $J_{\o_0}(\phi)=V^{-1}\sqrt{-1}\sum_{i=0}^{n-1}{n-i\over
n+1}\int_X\pl\phi\wedge\bar\pl\phi\wedge\o_\phi^i\wedge\o_0^{n-i-1}$,
$V=\int_X\o_0^n$, is the familiar Aubin-Yau functional.)
Returning to the problem at hand, we deduce that the $H_{(1)}(X\times
[0,1])$ Sobolev norms of the functions
$\phi(t;k)$ are uniformly bounded. The same is true for
the $H_{(1)}(X)$ Sobolev norms of $\phi(t;k)$ for each $t\in [0,1]$.

\bigskip
$\bullet$ The functions $\dot\phi(t;k)$ also satisfy an
interesting Harnack inequality of Li-Yau type.
Let $0\leq\tau<T\leq 1$ and let $\xi,X\in M$. Then we claim
\be
\dot\phi(\xi,\tau)\ \leq \dot\phi(X,T) + {1\over 8}\D(\xi,\tau,
X,T)
\ee
where
\be
\D(\xi,\tau, X,T)\ = \ \inf_\g \I_\tau^T\left({ds\over
dt}\right)^2\ dt
\ee
where $ds\over dt$ is the velocity in space at time $t$ and the
infimum is taken over all paths from $(\xi, \tau)$ to $(X,T)$
parametrized by $\tau\leq t\leq T$.
To see this,
let $L=2\dot\phi$. Then $\ddot\phi- |\pl\dot\phi|_{\o_\phi}^2 = |\pi_N
V|^2\geq 0$,
where $\pi_NV$ is the normal component of the holomorphic vector
field $V$ on ${\bf CP}^{N_k}$ generated by $\sigma^t$
\cite{PS02a}, and so
\be
{\pl L\over \pl t}\ \geq \ |D L|_\phi^2.
\ee
As in \cite{LY},
we can now choose a path $(t,s(t)) $ joining $(\xi,\tau) $  to  $(X,T) 
$ where
$\tau<T$ and $\xi,X\in M$. Then
\bea
L(X,T)-L(\xi,\tau)\ = \ \I_\tau^T {dL\over dt}\ dt\ &=& \
\I_\tau^T \left\{
{\pl L\over \pl t} + {\pl L\over \pl s}\cdot {ds\over
dt}\right\}\ dt\
\nonumber\\
&\geq&\ \I_\tau^T \left\{|DL|^2+{\pl L\over \pl s}\cdot
{ds\over
dt}\
\right\}\ dt.
\eea
Completing the square, we obtain
\be
L(X,T)-L(\xi,\tau)\ \geq \ -{1\over 4}\I \left({ds\over dt}
\right)^2\ dt.
\ee


\newpage

\begin{thebibliography}{99}

\bibitem{BT76} Bedford, E. and B.A. Taylor,
``The Dirichlet problem for a complex Monge-Ampre equation",
Invent. Math. {\bf 37} (1976), 1-44. 

\bibitem{BT82} Bedford, E. and B.A. Taylor,
``A new capacity for plurisubharmonic functions", Acta Math. {\bf 149} 
(1982), 1-40. 

\bibitem{B} Blocki, Z.,
``The complex Monge-Amp\`ere operator and pluripotential
theory", lecture notes available from the author's website.

\bibitem{C} Catlin, D.,
``The Bergman kernel and a theorem of Tian",
{\it Analysis and geometry in several complex variables} (Katata, 
1997), 1-23, 
Trends Math., Birkh\"auser Boston, Boston, MA, 1999. 

\bibitem{Ch} Chen, X.X.,
``The space of K\"ahler metrics", J. Differential Geom. 
{\bf 56} (2000), 189-234. 

\bibitem{CT} Chen, X.X. and G. Tian,
``Geometry of K\"ahler metrics and foliations by discs",
arXiv: math.DG / 0409433.

\bibitem{CLN} Chern, S.S., H. Levine,
and L. Nirenberg,
``Intrinsic norms on a complex manifold",
{\it Global Analysis, Papers in honor of K. Kodaira},
University of Tokyo Press (1969) 119-139.

\bibitem{D} Demailly, J.P.,
``{\it Complex analytic and differential geometry}",
book available from the author's website.

\bibitem{D97} Donaldson, S.K.,
``Remarks on gauge theory, complex geometry and $4$-manifold 
topology",
Fields Medallists' lectures, 384-403, World Sci. Ser. 20th Century 
Math., 
5, World Sci. Publishing, River Edge, NJ, 1997.

\bibitem{D99} Donaldson, S.K.,
``Symmetric spaces, K\"ahler geometry, and Hamiltonian
dynamics", Amer. Math. Soc. Transl. {\bf 196} (1999) 13-33.

\bibitem{D01} Donaldson, S.K.,
``Scalar curvature and projective imbeddings I",
J. Differential Geom. {\bf 59} (2001) 479-522.

\bibitem{D02} Donaldson, S.K.,
``Scalar curvature and stability of toric varieties",
J. Differential Geom. {\bf 62} (2002), 289-349.

\bibitem{D04} Donaldson, S.K.,
``Scalar curvature and projective imbeddings II",
arXiv: math.DG / 0407534.

\bibitem{GZ} Guedj, V. and A. Zeriahi,
``Intrinsic capacities on compact K\"ahler manifolds",
arXiv: math.CV / 0401302.

\bibitem{K} Klimek, M.,
{\it ``Pluripotential theory"}, London Mathematical Society
monographs, New Series {\bf 6} (1991) Oxford University Press,
New York.

\bibitem{L} Lelong, P.,
``Fonctions plurisousharmoniques et formes diff\'erentielles 
positives",
Gordon \& Breach, Paris-London-New York (1968).



\bibitem{LY} Li, P. and S.T. Yau,
``On the parabolic kernel of the Schr\"odinger operator", Acta Math. 
{\bf 156} (1986), 153-201.

\bibitem{Lu}  Lu, Zhiqin., ``On the lower order terms of the
asymptotic expansion of Tian-Yau-Zelditch",
Amer. J. Math.  122  {\bf 2} (2000),   235-273.


\bibitem{M87} Mabuchi, T.,
``Some symplectic geometry on compact K\"ahler manifolds",
Osaka J. Math. {\bf 24} (1987) 227-252.

\bibitem{P00} Paul, S.,
``Geometric analysis of Chow Mumford stability",
 Adv. Math. {\bf 182} (2004), 333-356.
  
\bibitem{PS02} Phong, D.H. and J. Sturm,
``Stability, energy functionals, and K\"ahler-Einstein
metrics", Comm. Anal. Geometry {\bf 11} (2003)
563-597, arXiv: math.DG / 0203254.

\bibitem{PS02a} Phong, D.H. and J. Sturm,
``Scalar curvature, moment maps, and the Deligne pairing",
Amer. J. Math. {\bf 126} (2004) 693-712,
arXiv: math.DG / 
0209098.

\bibitem{PS03} Phong, D.H. and J. Sturm,
``The Futaki invariant and the Mabuchi energy of a complete
intersection", Comm. Anal. Geometry {\bf 12} (2004)
321-343, arXiv: math.DG / 0312529.

\bibitem{S92} Semmes, S.,
``Complex Monge-Amp\`ere and symplectic manifolds",
Amer. J. Math. {\bf 114} (1992) 495-550.

\bibitem{T90} Tian, G.,
`` On a set of polarized K\"hler metrics on algebraic
manifolds", J. Diff. Geom. {\bf 32} (1990) 99-130.


\bibitem{T97} Tian, G.,
``K\"ahler-Einstein metrics with positive scalar curvature",
Invent. Math. {\bf 130} (1997) 1-37.

\bibitem{Y78} Yau, S.T.,
``On the Ricci curvature of a compact K\"ahler manifold
and the complex Monge-Amp\`ere equation I",
Comm. Pure Appl. Math. {\bf 31} (1978) 339-411.

\bibitem{Y87} Yau, S.T.,
`` Nonlinear analysis in geometry",
Enseign. Math. (2)  {\bf 33}  (1987),  no. 1-2, 109--158.


\bibitem{Y} Yau, S.T.,
``Open problems in geometry",
Proc. Symp. Pure Math. {\bf 54}, AMS Providence, RI
(1993) 1-28.

\bibitem{Z} Zelditch, S.,
``The Szeg\"o kernel and a theorem of Tian",
Int. Math. Res. Notices {\bf 6} (1998) 317-331.

\bibitem{Z96} Zhang, S.,
``Heights and reductions of semi-stable varieties",
Compositio Math. {\bf 104} (1996) 77-105.


\end{thebibliography}
\enddocument